\documentclass{birkart}

 \newtheorem{thm}{Theorem}[section]

 \newtheorem{prop}[thm]{Proposition}
 \theoremstyle{definition}

 \theoremstyle{remark}
 \newtheorem{rem}[thm]{Remark}

\numberwithin{equation}{section}
\numberwithin{figure}{section}

\newcommand{\C}{{\mathbb C}}
\newcommand{\D}{{\mathbb D}}

\newcommand{\R}{{\mathbb R}}
\newcommand{\s}{{\mathbb S}}
\newcommand{\Z}{{\mathbb Z}}

\newcommand{\schlicht}{{\mathcal S}}

\newcommand{\diff}{{\mathrm d}}

\newcommand{\tor}{{\boldsymbol\Pi}}
\newcommand{\Delt}{{\mathbf D}}
\newcommand{\sur}{{\mathbf S}}
\newcommand{\bp}{{\mathbf p}}
\newcommand{\bq}{{\mathbf q}}
\newcommand{\proj}{\boldsymbol \pi}
\newcommand{\btau}{\boldsymbol \tau}
\newcommand{\bsigma}{\boldsymbol \sigma}
\newcommand{\bo}{\boldsymbol\Omega}
\newcommand{\binfty}{\boldsymbol\infty}
\newcommand{\bzero}{\boldsymbol 0}
\newcommand{\bvarphi}{\boldsymbol\varphi}
\newcommand{\bphi}{\boldsymbol\phi}
\newcommand{\bpi}{\boldsymbol{\mathcal D}}
\newcommand{\bgamma}{\boldsymbol\Gamma}
\newcommand{\iu}{{\rm i}}
\newcommand{\bbpsi}{\boldsymbol\Psi}
\newcommand{\funddomain}{{\mathcal D}}

\begin{document}
%
\title[Branch point area methods]
{Branch point area methods in \\conformal mapping}
\author[Abuzyarova]
{Natalia Abuzyarova}

\address{Abuzyarova: Department of Mathematics\\
The Royal Institute of Technology\\
S -- 100 44 Stockholm\\
Sweden}

\email{naab@math.kth.se}

\thanks{Both authors wish to thank the G\"oran Gustafsson 
Foundation for generous support.}

\author[Hedenmalm]{H\aa{}kan Hedenmalm}
\address{Hedenmalm: Department of Mathematics\\
The Royal Institute of Technology\\
S -- 100 44 Stockholm\\
Sweden}
\email{haakanh@math.kth.se}

\subjclass{Primary 99Z99; Secondary 00A00}

\keywords{Class file, Birkart}


\begin{abstract}
The classical estimate of Bieberbach -- that $|a_2|\le2$
for a given univalent function $\varphi(z)=z+a_2z^2+\ldots$ in the class
$\schlicht$ -- leads to best possible pointwise estimates of the ratio 
$\varphi''(z)/\varphi'(z)$ for $\varphi\in\schlicht$, first obtained by
K\oe{}be and Bieberbach.
For the corresponding class $\Sigma$ of univalent functions in the
exterior disk, Goluzin found in 1943 -- by extremality methods -- the 
corresponding  best possible pointwise estimates of $\psi''(z)/\psi'(z)$ 
for $\psi\in\Sigma$. 
It was perhaps surprising that this time, the expressions involve
elliptic integrals. Here, we obtain the area-type theorem which has
Goluzin's pointwise estimate as a corollary. This shows that the
K\oe{}be-Bieberbach estimate as well as that of Goluzin are both firmly 
rooted in the area-based methods. The appearance of elliptic integrals finds 
a natural explanation: they arise because a certain associated covering 
surface of the Riemann sphere is a torus.
\end{abstract}

\maketitle

\section{Introduction}

{\bf Area methods.} Area methods play an important role in the theory of 
conformal mappings.
The original Gr\"onwall area theorem states that if $\psi$ belongs to
the class $\Sigma$, with series expansion
$$\psi(z)=z+\sum_{n=0}^{+\infty}b_n\,z^{-n},$$
then
\begin{equation}
\frac1{\pi}\int_{\D_e}\big|\psi'(z)-1\big|^2\,\diff A(z)
=\sum_{n=0}^{+\infty}n\,|b_n|^2\le1.
\label{eq-areathm}
\end{equation}
Here, $\diff A(z)=\diff x\diff y$ is ordinary area measure in the plane.
Also, we recall that $\psi\in \Sigma$ means that $\psi$ is a conformal
mapping from the exterior disk
$$\D_e=\big\{z\in\C\cup\{\infty\}:\,1<|z|\le+\infty\big\}$$
to some domain on the Riemann sphere $\s =\C_{\infty}$, with the properties
that $\psi(\infty)=\infty$, and $\psi'(\infty)=1$. In particular,
(\ref{eq-areathm}) implies that $|b_1|\le1$. After an inversion
of the plane plus a square root transformation, it follows that for
$\varphi$ in the class $\schlicht$ of conformal mappings of the unit disk
$\D$ into $\C$ with $\varphi(0)=0$ and $\varphi'(0)=1$, we have the
estimate $|\varphi''(0)|\le4$. The M\oe{}bius automorphisms of the unit disk
allow us to move the point at the origin to an arbitrary point in $\D$; this
results in the K\oe{}be-Bieberbach estimate
\begin{equation}
\left|\frac{\varphi''(z)}{\varphi'(z)}-\frac{2\bar z}{1-|z|^2}\right|
\le\frac{4}{1-|z|^2},\qquad z\in\D.
\label{eq-ptest}
\end{equation}
This estimate is best possible in the sense that if we consider, for a given
$z_0\in\D$, the set of points
$$\left\{\frac{\varphi''(z_0)}{\varphi'(z_0)}:
\,\,\varphi\in\schlicht\right\},$$
we obtain a closed circular disk of radius $4/(1-|z_0|^2)$ centered at
$2\bar z_0/(1-|z|^2)$.
\bigskip

\noindent
{\bf Goluzin's inequality.} For the class $\Sigma$, Goluzin 
\cite{Gol1943}, \cite[p. 132]{Goluz} found in 1943 the estimate analogous 
to (\ref{eq-ptest}) using extremality methods. Given $\psi\in\Sigma$, 
it reads:
\begin{equation}
\left|\frac{\psi''(z)}{\psi'(z)}+\frac{4|z|^2-2}{z(|z|^2-1)}
-\frac{4\bar z}{|z|^2-1}\,\frac{E\left(\frac1{|z|}\right)}
{K\left(\frac1{|z|}\right)}\right|
\le\frac{4|z|}{|z|^2-1}\,\left(1-\frac{E\left(\frac1{|z|}\right)}
{K\left(\frac1{|z|}\right)}\right),
\label{eq-ptests}
\end{equation}
for $z\in\D_e$. Here, $E$ and $K$ are the elliptic integrals
\begin{equation}
E(\lambda)=\int_0^1\sqrt{\frac{1-\lambda^2t^2}{1-t^2}}\,\diff t,\qquad
\lambda\in\D,
\label{el-E}
\end{equation}
and
\begin{equation}
K(\lambda)=\int_0^1\frac{\diff t}{\sqrt{(1-\lambda^2t^2)(1-t^2)}},
\qquad\lambda\in\D.
\label{el-K}
\end{equation}
Like (\ref{eq-ptest}), the estimate (\ref{eq-ptests}) is best possible.
However, the derivation of (\ref{eq-ptests}) which Goluzin employs is
quite different from the above-mentioned classical derivation of
(\ref{eq-ptest}) in terms of area estimates. Here, {\sl we find the area-type
estimate needed to derive} (\ref{eq-ptests}). Basically, we introduce
a square root slit in $\s$ between the point at infinity and a given point
$\psi(z_0)$ for $z_0\in\D_e$, and apply Stokes' theorem to the resulting
compact covering surface over the Riemann sphere. The application of Stokes'
theorem involves the use of the Green function for the part of the covering
surface which covers $\psi(\D_e)$; in terms of the coordinates of $\D_e$, this
Green function results from applying a square root slit in $\D_e$ between
infinity and $z_0$. This latter surface is conformally equivalent to an
annulus.
From the area-type method point of view which is hinted above and described
in detail in the following sections, the Green function for the annulus --
which is expressible in terms of elliptic integrals -- is the reason why
elliptic integrals appear in (\ref{eq-ptests}).

In the paper of Bergman and Schiffer \cite{BerSch}, the reader will find
out how the Grunsky inequalities (a general version of the area method) 
can be developed from the perspective of Bergman kernels. Also, he (or she)
may find it interesting to compare the branch point area methods that are
developed here with those of Nehari \cite{Neh2}. 

\section{The area-type inequality}

{\bf An application of Stokes' theorem.}
Let $\sur$ be a compact Riemann surface. We will later consider the special 
case when $\sur$ is a
(branched) covering surface of the Riemann sphere $\s =\C_{\infty }$.
The Sobolev space $W^{1,2} (\sur)$ consists of those locally summable 
functions $f:\ \sur\to\C$
for which  the first-order differential $\omega_f =\diff  f$
is an element of the Hilbert space of $1$-forms $L^2_1 (\sur)$ (see 
\cite[Ch. 7, pp. 181--182]{Springer}). We recall the standard definition
of the norm in $L^2_1 (\sur )$:
$$
\| \omega \|^2_{L^2} =\int_{\sur} \omega \wedge {}^*\bar\omega .
$$
Here, we use the standard Hodge notation
$$
\omega=u\,\diff z+v\,\diff\bar z,\qquad
{}^*\omega=-\iu u\, \diff z+ \iu v\, \diff\bar z,
$$
where $z$ is any local complex parameter.
The space $W^{1,2} (\sur)$ is supplied with the semi-norm
$$
\| f\|^2_{W^{1,2}} =\| \diff f\|^2_{L^2} .
$$
We will consider the space $W^{1,2}(\sur)$ as taken modulo the constant
functions; that is, any constant function will be thought of as the zero
function. This is done with the intention to make the above semi-norm a 
norm on $W^{1,2}(\sur)$.
In terms of a local complex parameter $z$, the differential 
$\omega_f =\diff f$ may be written as
$$\omega_f =\partial_z f\,\diff z+\bar{\partial}_z f\,\diff\bar z.$$
This is the local form of the global decomposition
$$
\omega_f =\omega_{f,1} +\omega_{f,2},
$$
where in terms of local coordinates  $\omega_{f,1} =\partial_z f\,\diff z$,   
$\omega_{f,2} =\bar{\partial}_z f\,\diff \bar z$
(see \cite[Ch. 1, p. 62--63 and Ch. 2, p. 153]{Forst}).

The function $f\in W^{1,2} (\sur)$  generates the second-order differentials
$$\Lambda_{f,1} = \omega_{f,1}\wedge\bar\omega_{f, 1} , \quad
\Lambda_{f,2} =-\omega_{f,2}\wedge \bar\omega_{f,2},$$
which have the form
\begin{equation}
\Lambda_{f,1} =|\partial_z f |^2\, \diff z\wedge \, \diff\bar{z},
\qquad \Lambda_{f,2} =|\bar\partial_{z} f |^2\, \diff z\wedge
\, \diff \bar{z} ,
\label{2-0}
\end{equation}
in a local complex parameter $z$.
Note that
$$
\| f \|^2_{W^{1,2}} =\iu\int\limits_{\sur} \Lambda_{f,1} +
\iu\int\limits_{\sur} \Lambda_{f,2} .
$$

The next result is a consequence of Stokes' theorem.

\begin{prop}
\it For $f\in W^{1,2} (\sur),$ both integrals
$\int\limits_{\sur}\Lambda_{f,1}$
and $\int\limits_{\sur}\Lambda_{f,2}$ are finite, and
\begin{equation}
\int\limits_{\sur}\Lambda_{f,1}
=\int\limits_{\sur}\Lambda_{f,2}.
\label{2-1}
\end{equation}
\label{prop2-1}
\end{prop}

\begin{proof}
Assume that $f\in C^2 (\sur ),$ and consider the integral
$$
\int\limits_{\sur} \diff\bigl( f\, \diff\bar{f}\bigr) .
$$
Simple calculations give us
$$
\diff \bigl(f\, \diff\bar{f}\bigr) =
\bigl(|\partial_z f|^2-|\bar\partial_z f|^2\bigr) \, \diff z\wedge
\, \diff\bar{z}
$$
in a local complex parameter $z$. This means that
\begin{equation}
\int\limits_{\sur} \diff \bigl(f\, \diff\bar{f}\bigr) =
\int\limits_{\sur} \Lambda_{f,1}
-\int\limits_{\sur}\Lambda_{f,2} .
\label{2-2}
\end{equation}
By Theorem 6-4 \cite[Ch. 6, p. 167]{Springer}, we have
$$
\int\limits_{\sur}
\diff\bigl( f\, \diff\bar{f} \bigr) =0 .
$$
In view of (\ref{2-2}), we obtain
$$
\int\limits_{\sur} \Lambda_{f,1} =
\int\limits_{\sur} \Lambda_{f,2}.
$$

The general case $f\in W^{1,2}(\sur )$ follows by approximation argument.
\end{proof}

We point out that Proposition~\ref{prop2-1} claims the following:
for the exact first-order differential form $\omega=\omega_f$,
$$
\int\limits_{\sur} \omega\wedge \bar{\omega} =0,
$$
which, of course, is not true for an arbitrary $1$-form.

\bigskip

\noindent
{\bf Solution of Laplace's equation on a subdomain.}
We consider a nontrivial finitely connected subdomain $\bo$ of the compact
Riemann surface $\sur$ (nontriviality means that $\bo\neq\emptyset,\sur$), 
and a meromorphic  function $R$ on $\sur$, the poles of which are all 
contained in $\bo$. The poles of $R$ are denoted by $p_1,\dots, p_N$, and 
$m_j$ is the order of the pole $p_j$, for $j=1,\ldots,N$.

\begin{prop}
There exists a function $Q:\, \sur\to\s$ with the following properties:

\noindent
{\rm(Q1)} $Q$ equals zero on $\sur\setminus\bo ;$

\noindent
{\rm(Q2)} $Q$ is harmonic on $\bo\setminus\{ p_1,\dots ,p_N\}; $

\noindent
{\rm(Q3)} the function $P=R-Q$ is of H\"older class $\text{\rm Lip}\,\frac12$
on $\sur$, and it belongs to the Sobolev space $W^{1,2} (\sur).$
\label{prop2-3}
\end{prop}

\begin{proof} As a matter of convenience, we assume in the first part of
this proof that the domain $\bo$ has real-analytic boundary.
For  $\bo$, considered as a Riemann surface, we introduce its {\it conjugate 
surface} $\bo^{*} $  (see \cite[Ch. 8, p. 217, Problem 1]{Springer}). Let 
$\bo^{*}$ be another copy of $\bo$
and $*:\ \bo\to\bo^{*} $ be the identity mapping, $p^*=*(p)$. We also use
the same notation $*$ for the inverse mapping, $*=*^{-1}$, so that $p^{**}
=p$. The complex structures of $\bo$ and $\bo^*$ are different, however: 
if $z=\Phi (p)$ is a local complex parameter about some point $p_0\in\bo$,
with $\Phi(p_0)=0$, we pick $\bar z=\bar\Phi (p)=\Phi^*(p^{*})$ as
a local complex parameter about $p_0^{*}$, where the latter relation is 
used to define the function $\Phi^*$. Out of $\bo$ and $\bo^{*}$, we form 
the  {\sl Schottky double}
$$\widehat{\bo} =\bo\cup\bo^*\cup\partial\bo$$
by identifying conjugate boundary points $p\in\partial \bo$ and 
$p^{*}\in\partial \bo^{*}.$
As a local complex parameter near the identified boundary points 
$p_0=p_0^{*}\in \partial\bo$, we pick
$$
z=\begin{cases}
\Phi (q),\quad p\in \bo\cup\partial\bo,\\
\bar\Phi(p^*) ,\quad  p\in\bo^{*},
\end{cases}
$$
where $z=\Phi(p)$ is a special type of local complex parameter about $p_0$:
it is defined on some neighborhood $V\subset \sur$ around $p_0$, and it 
maps $V\cap\bo$ onto a region in the upper half-plane $\rm{Im}\,z>0$, with 
$\Phi (p_0) =0$ and such that the connected segment of $\partial \bo\cap V$ 
containing $p_0$ is mapped onto a segment of the real axis 
(see \cite[Ch. 8, p. 217, Problem 2]{Springer}). This way, we supply 
$\widehat\bo$ with the structure of a compact Riemann surface. 
By Corollary 8-1 in \cite[Ch. 8, p. 211]{Springer}, for every 
point $p_j$, there exist functions $g_j$ and $g_j^{*}$  with the following
properties:

\noindent
$\bullet$ $g_j$  is harmonic in $\widehat{\bo}\setminus \{p_j\}$, and
$g_j^*$  is harmonic in $\widehat{\bo}\setminus \{p_j^*\}$;

\medskip

\noindent
$\bullet$ $g_j$ has at the point $p_j$ the same singularity as $R$, while
$g_j^*$ has at the point $p_j^* $ the same singularity as $-R\circ *$.

\medskip

\noindent We now put 
$$
Q_j(p)=\frac{1}{2}\, \Big\{ g_j(p)+g_j^{*}(p)
-g_j(p^{*})-g_j^{*}(p^{*})\Big\} .
$$
The function $Q_j$ has the following properties, for $j=1,\ldots,N$:

\medskip

\noindent
(1) it is harmonic in $\bo\setminus\{ p_j\} ;$

\medskip

\noindent
(2) the function $R-Q_j$ is regular at the point $p_j;$

\medskip

\noindent
(3) $Q_j$ is continuous in $\bar{\bo}\setminus\{p_j\}$,
and $Q_j(p)=0$ for $p\in\partial \bo$.

\bigskip

\noindent Next, we define the function $Q$ by
$$
Q (p) =
\begin{cases}
\sum\limits_{j=1}^{N} Q_j (p), \qquad   p\in\bo ,\\
    0,\qquad  p\in \sur\setminus \bo ,
\end{cases}
$$
and introduce the associated function $P$, as given by
$$
P(p)=R (p)-Q (p).
$$
The properties of $Q$ imply that $P$ coincides with $R$ on the compact set 
$\sur\setminus \bo$, and that $P$ extends harmonically across the set 
$\{ p_1,\dots ,p_N\}$. Moreover, in view of the real-analyticity of the 
boundary $\partial \bo$ it follows that the function $Q$ is 
Lipschitz-continuous near $\partial\bo$, making $P$ Lipschitz-continuous 
on all of $\sur$, and hence we get $P\in W^{1,2}(\sur)$.

All the above considerations are valid under the assumption that $\bo$ has
real-analytic boundary. In the general case, we may approximate $\bo$ by a 
increasing sequence of domains $\bo_n$ with real-analytic boundaries.
For each such domain $\bo_n$, we construct the function $Q_n$ according to
the above scheme. We then appeal to a well-known result of Beurling 
\cite[p. 53]{Ber}, which implies the uniform boundedness of the local  
Lip$\,\frac{1}{2}$-norms (away from the poles of $R$ $\{ p_1, \dots , p_N\}$) 
of $Q_n.$ Thus, the sequence $\{Q_n \}$ converges in a weak sense to some 
function $Q$, defined on $\bo$. We set $P=R-Q$  with this limit function 
$Q$.
The functions $P$ and $Q$ satisfy all required conditions, with one possible 
exception: we need to show that $P\in W^{1,2} (\sur ).$
However, this is an obvious consequence of the following fact: the function 
$P$ solves the Dirichlet problem on $\bo$ with boundary values equal to $R$, 
and the solution to the Dirichlet problem minimizes the Dirichlet integral 
over $\bo$. The $W^{1,2} (\sur )$-(semi-)norm of $P$ is the sum of its 
Dirichlet integral over $\bo$ and the Dirichlet integral of $R$ over 
$\sur\setminus\bo $, which both are finite.
In view of this, we conclude that $P$ belongs to $ W^{1,2} (\sur )$.
\end{proof}

\medskip

\noindent
{\bf The area-theorem type inequality.}
We want to apply (\ref{2-1}) to $P=R-Q.$ For this function, we have, by 
(\ref{2-0}),
$$
\Lambda_{P,1} =\big|\partial_z (R-Q) \big|^2
\, \diff z\wedge \, \diff\bar{z} ,\quad
\Lambda_{P,2} = \big|\bar\partial_z(R-Q)\big|^2 \, 
\diff z\wedge \, \diff \bar{z},
$$
where $z$ is a local complex parameter.

Note that the area element $ \diff A (z)$ is $\frac{\iu}{2}\, 
\diff z\wedge \, \diff\bar{z}$. We have
$$
\frac{\iu}{2} \, \int\limits_{\sur} \Lambda_{P,1}\ge
\frac{\iu}{2}\, \int\limits_{\bo} \Lambda_{P,1}
$$
and
$$
\int\limits_{\sur} \Lambda_{P,2} =
\int\limits_{\bo} \Lambda_{P,2}.
$$
Combining these relations with (\ref{2-1}), applied to the function $P,$ 
we obtain
\begin{equation}
\int\limits_{\bo} \frac{\iu}{2}\Lambda_{P,1} \le
\int\limits_{\bo} \frac{\iu}{2}\Lambda_{P,2},
\label{2-3}
\end{equation}
where, in terms of local coordinates,
$$
\frac{\iu}{2}\, \Lambda_{P,1} = \big|\partial_z R -\partial_zQ \big|^2 \, 
\diff A (z) ,
\quad
\frac{\iu}{2}\, \Lambda_{P,2} =\big|\bar\partial_zQ \big|^2\, \diff A (z) .
$$

Let us note that we have equality in (\ref{2-3}) precisely when the 
complement $\sur\setminus\bo$ has zero area.

In the next section we will consider more concrete choice of  $\sur,$ $\bo$
and $P,$ to derive from (\ref{2-3}) area theorem type estimates for 
univalent functions.

\section{Applications}

{\bf The torus subdomain.}
For a covering surface $\sur$ of the Riemann sphere, we will denote by 
$\proj$ the projection mapping of $\sur$ onto $\s$.

Let $\sur$ be the image of the Riemann sphere $\s$
under the mapping $z\mapsto z^2.$
Thought of as a covering surface of $\s,$ $\sur$ is a two-sheeted covering, 
with associated projection $\proj :\sur\to\s .$
The covering  has two branch points in $\sur ,$ which we call 
$\boldsymbol 0$ and $\boldsymbol\infty.$
They project to the points $0$ and $\infty :$ $\proj (\boldsymbol 0) =0$ 
and $\proj (\boldsymbol \infty) =\infty .$

We now describe a concrete domain $\bo .$ Let $\varphi (w)$ be a univalent 
function, defined in the unit disk $\D$, which maps into $\s,$
such that for some real parameter $x_0$, $0<x_0<1,$ we have
$$\varphi (x_0)=0,\quad  \varphi (-x_0 )=\infty ,\quad
\varphi' (x_0 )=1.$$
We put $\Omega=\varphi (\D )$ and  note that $\Omega $ contains the points 
$0$ and $\infty.$ We will use the notation $\phi $ for the inverse function 
to $\varphi :$
$$
\phi=\varphi^{-1} :\, \Omega\to\D .
$$

Denote by $\bo$ the lifting of  $\Omega$ to $\sur,$ so that  
$\proj (\bo )=\Omega.$ To get $\bo $, we  should first cut $\Omega$ from 
$0$ to $\infty,$  then take two copies of such cut $\Omega ,$ and attach 
them crosswise along the cuts.

The preimage of the cut from $0$ to  $\infty$  in $\Omega$ is a cut from 
$x_0$ to $-x_0$ in the unit disk $\D$. Attaching crosswise along these 
preimage cuts two replicas of cut $\D$, we get a two-sheeted covering 
surface $\Delt$, which is conformally equivalent to  $\bo$.
The surface $\Delt$  has two branch points, which  project to the 
points $x_0$ and $-x_0$ of the unit disk.

We need to define an analytic self-mapping $\Delt\to\Delt$. It will be the
correspondence $p\mapsto p'$ between the points $p$ and $p'$ belonging to 
the different sheets of $\Delt$. Namely, the point $p$ with the projection 
$\proj (p)=z$, $z\in \D\setminus \{ x_0, \ -x_0\}$, is mapped to another 
point $p'\in \Delt$  with the same projection  $\proj (p')=z$.
For $p$ such that $\proj (p) =\pm x_0 ,$ we put $p'=p$. We will call $p'$ 
{\sl the mirror point} to the point $p$.

We now define the mapping $\bvarphi :\, \Delt\to\bo$
to be the lifting of $\varphi$ to $\Delt$. By definition, it maps the 
point $p\in\Delt$ to the point $\bp\in\bo$ with the projection 
$\varphi (\proj (p)),$ $i=1,\, 2.$
We also define
$$\bphi  =\bvarphi ^{-1} :\, \bo\to\Delt,$$
which is the lifting of $\phi:\Omega\to\D$.

Further, denote by $\bp ' $ the point $\bvarphi (p')$ for 
$\bp =\bvarphi (p)\in\bo .$  We get an analytic self-mapping
$\bo\to\bo$, which takes any point $\bp\in\bo\setminus 
\{ \bzero ,\binfty\}$ to the other point
$\bp'\in\bo\setminus\{\bzero ,\binfty\}$ with the same projection: 
$\proj (\bp )=\proj (\bp' );$
 each of the points $\bzero ,$ $\binfty ,$ is taken to itself. As in the 
case of the points $p,\, p'\in\Delt$, we will call $\bp'$ {\sl the mirror 
point} to the point $\bp$.

Next, we introduce a meromorphic function $R(\bp ),$ $\bp\in \sur,$ which 
has a simple  pole at the branch point $\bzero$ and has no other poles. 
Note that any meromorphic function $f$ on our surface $\sur $ can be 
expressed in terms of the global coordinates of $\s =\C\cup\{\infty\}$ as
$$
f(z)=f_1(z)+\sqrt{z} f_2(z),
$$
where $f_1$ and $f_2$ are  meromorphic functions on $\s$,
and $\sqrt{z}$ means the algebraic square root of $z$.
{\sl We define the function $R$ to be the above  $f$ with the choices 
$f_1(z)=0$, $f_2(z)=1/z$.}

Our next project is to construct the function  $Q$, which
satisfies the conditions (Q1)--(Q3) of Proposition~\ref{prop2-3} for
this given $R$. To this end, as a first step, we consider the Green function
$G_{\bo} (\bp ,\bq)$ of the domain $\bo$. 
For fixed  $\bq\in\bo$, the  function $\bp\mapsto G_{\bo} (\bp ,\bq )$ 
is harmonic on $\bo\setminus\{ \bq \}$, vanishes on the boundary 
$\partial\bo ,$ and has the logarithmic singularity $-\log |z| + O(1)$ in 
terms of local coordinates around  $\bp =\bq$.
The function $\bphi$ maps $\bo$ onto $\Delt$ conformally.
It follows that (see \cite[Ch. 6, \S 2, pp. 201--202]{Nevanl})
$$
G_{\bo} (\bp ,\bq ) =G_{\Delt} (\bphi (p),
\bphi (q)) ,  \qquad p,q\in\Delt ,\quad \bp =\bphi (p), \quad \bq =\bphi (q).
$$

For $p,q\in\Delt$,  we define
\begin{eqnarray*}
G_{\Delt}^{\rm alt} (p,q) &=&G_{\Delt}(p,q)-G_{\Delt} (p',q) , \\
G_{\bo}^{\rm alt} (\bp ,\bq )&=& 
G_{\Delt}^{\rm alt} (\bphi(p ),\bphi(q )), \qquad \qquad  
\bp =\bphi (p), \,\,\, \bq =\bphi (q).
\end{eqnarray*}
From the above definitions, it follows that
\begin{eqnarray}
G_{\bo}^{\rm alt} (\bp ',\bq ) &=&-G_{\bo}^{\rm alt} (\bp ,\bq ),
\qquad \bp,\, \bq \in\bo ,\nonumber \\
G_{\Delt}^{\rm alt} (p',q) &=&-G_{\Delt}^{\rm alt} (p,q),\qquad\,\,
p,\,q\in\Delt.
\label{3-2}
\end{eqnarray}

The functions $R$ and $R_{\Delt} =R\circ \bvarphi$ have
the same property:
\begin{eqnarray*}
R (\bp ')&=&-R (\bp ),\qquad \bp ,\ \bp '\in \sur, \\
R _{\Delt}(p')&=&-R _{\Delt}(p),\qquad p,\ p'\in \Delt .
\end{eqnarray*}

\medskip

For our further considerations, we need yet another covering surface of 
$\s$. To get it, we first supply  $\s$ with two cuts. One of the cuts is 
made in $\D$ from $-x_0$ to $x_0$ as we did earlier for the description of  
$\Delt$.
The second cut goes from $-1/x_0$ to $1/x_0$. Also, this cut is to be 
obtained from the first one by reflection in the unit circle: 
$z\mapsto 1/\bar{z}$.
As we attach two copies of such cut Riemann spheres  crosswise along the 
corresponding (same) cuts, we obtain a compact surface, which we denote by 
$\tor$. It is a two-sheeted covering surface of $\s$ with four branch points.
In terms of conformal equivalence, $\tor$ is a torus.
As the second cut from $-1/x_0 $ to $1/x_0$ falls outside the unit disk 
$\D$, we may think of the surface $\Delt$ as a subdomain of  $\tor $.

For a moment, let us fix an arbitrary $q\in\Delt.$ In addition to (\ref{3-2}),
$G_{\Delt}^{\rm alt}(p,q)$ has the following properties:

\medskip

\noindent
(1) $G_{\Delt}^{\rm alt}(p,q) =0,$ for $ p\in\partial \Delt .$

\medskip

\noindent
(2) the function $p\mapsto G_{\Delt}^{\rm alt} (p,q)$ has the logarithmic 
singularity $-\log |z|+O(1)$ in terms of local coordinates around  
the point $p=q$
and the logarithmic singularity $\log |z|+O(1)$ in terms of local 
coordinates around  the point $p=q'$, where $q,$ $q'$  are mirror points 
to each other (so that $q\neq q'$ and $\proj (q) =\proj (q')\in\D$).

\medskip

\noindent
(3) $G_{\Delt}^{\rm alt}(p,q)$ is harmonic on $\bar{\Delt} \setminus 
\{ q\}$.

\bigskip

Next, we describe a self-mapping $\tor\to\tor$, {\sl reflection in } 
$\partial\Delt$.
Namely, this mapping takes the point $p$ with the projection $\proj (p) =z$ 
to the point $p^*$ with the projection
$\proj (p^*) =1/\bar z.$ The choice of $p^*$ from the two different points 
of $\tor$ with the same projection $1/\bar z$
is defined by the following requirements: our mapping must be continuous on 
$\tor$, and $p^*=p$ for all $p\in\partial \Delt$.

The function $G_{\Delt}^{\rm alt}$ may be extended harmonically across the 
boundary $\partial\Delt$. Indeed, by the Schwarz reflection principle, for 
any fixed $q\in\Delt$, we define $G_{\Delt}^{\rm alt}(p,q)$ on the 
complement of $\Delt$ by
$$
G_{\Delt}^{\rm alt} (p,q) =-G_{\Delt}^{\rm alt} (p^*,q), 
\qquad p\in\tor\setminus\Delt\setminus\{ q^*,\, (q')^{*}\} .
$$
The extended function $p\mapsto G_{\Delt}^{\rm alt}(p,q)$ is harmonic on 
$\tor  \setminus\{q,\ q',\ q^{*},\ (q')^{*}\}.$ It has the singularity 
$-\log |z|+O(1)$ in terms of local coordinates around the points 
$p=q$ and $p=(q')^*$,
and the singularity $\log |z|+O(1) $ in terms of local coordinates around  
the points $p=q'$ and $p=q^{*}$.

\begin{rem}
\rm {The reason why we consider the Green functions $G_{\bo}$, $G_{\Delt}$ 
and the functions $G_{\bo}^{\rm alt},$ $G_{\Delt}^{\rm alt} ,$ is given by 
the following observation. Let $\Omega$ be a subdomain of $\s$ with analytic 
boundary. Also, we assume that $\Omega $ contains $0$ and has the property: 
$w\in \Omega\Longleftrightarrow -w\in\Omega$. We introduce
$$
G_{\Omega}^{\rm alt } (w,\lambda )=G_{\Omega} (w,\lambda )-
G_{\Omega} (-w,\lambda ),
$$
where $G_{\Omega}$ is the Green function of $\Omega$. This function can 
be represented as
$$
G_{\Omega}^{\rm alt} (z,\lambda ) =-\frac{1}{2}\log |w-\lambda |^2
+\frac{1}{2}\log|w+\lambda |^2 +H(w,\lambda ),
$$
where $H(w,\lambda )$ is an odd harmonic function of the variable $w$.
We observe that   the function $Q$ defined by
$$
Q(w)=\partial_{\lambda} G_{\Omega}^{\rm alt} (w,\lambda )\vert_{\lambda =0} =
\frac{1}{w} +\partial_{\lambda} H(w,\lambda )\vert_{\lambda =0}
$$
is harmonic on $\bar{\Omega}\setminus\{ 0 \}$, with a simple pole at 
the point $w=0$, and it equals zero on $\partial \Omega .$}
So, in the special case $\sur =\s,$ $R(z)=1/z,$ we obtain the function 
$Q$ of Proposition~\ref{prop2-3} from the function $G_{\Omega}^{\rm alt}$ 
in the above manner.
\end{rem}

We shall try to find the required  function $Q:\, \sur\to\s$ for the given 
$R:\, \sur\to\s$ in an analogous fashion.
The theory of elliptic functions (or integrals) is needed to obtain the 
explicit form of $G^{\rm alt}_{\Delt}.$

\bigskip

\noindent
{\bf Elliptic functions and the Green function for the torus subdomain.}
We recall some definitions and facts from the elliptic functions theory  
(see \cite[Ch.V, VI]{Akhiez}).

Let $k$ be a real parameter, $0<k<1$. We introduce the following notation:
\begin{equation}
\begin{array}{ll}
{k}'=\sqrt{1-k^2} ,\qquad l=\frac{1-k'}{1+k'} ,\qquad l'=\sqrt{1-l^2} , 
\qquad M=\frac{1}{1+k'},
\\
K=K(k), \qquad K'=K(k'), \qquad L=K(l) ,\qquad L'=K(l') ,
\end{array}
\label{M}\end{equation}
where the function  $K(\lambda )$ is defined by  (\ref{el-K}).
The values $L,$ $L',$ $K,$  $K'$ are connected by Landen's transformation 
(see \cite[Ch.VI]{Akhiez}), namely,
\begin{equation}
K=2ML,\quad K'=ML' .
\label{L-M}
\end{equation}

Let $h=e^{-\pi K'/K}.$  One of Jacobi's theta-functions
$\vartheta_0(u)$ is defined by
$$
\vartheta_0(u) =1-2h\cos 2\pi u+2h^4\cos 4\pi u-2h^9\cos 6\pi u+\dots ,
\qquad u\in\C .
$$
We also recall the definitions of the following Jacobi elliptic functions:
\begin{gather}
\theta_0 (z) =\vartheta_0\left( \frac{z}{2K}\right), \qquad
Z(z)=\frac{\theta_0'(z)}{\theta_0 (z)},
\label{el-J}
\\
{\rm sn}  (z;k) =\frac{\iu e^{-\frac{\pi \iu}{4K}(2z+\iu K')}}{\sqrt{k'}}\, 
\frac{\theta_0 (z-\iu K')}{\theta_0 (z)} ,
\qquad {\rm dn}  (z;k) =\sqrt{k'}\, \, \frac{\theta_0 (z-K)}{\theta_0 (z)},
\nonumber
\\ {\rm cn}  (z;k) =-\iu e^{-\frac{\pi \iu}{4K}(2z+\iu K')}\, 
\sqrt{\frac{k'}{k}}\, \frac{\theta_0 (z-K-\iu K')}{\theta_0 (z)} .
\nonumber
\end{gather}

Let $\kappa =\frac{2x_0}{1+x_0^2}$ for some real $x_0,$ $0<x_0<1 .$ The 
point $x_0$ is the same one we used to define the surfaces $\Delt$,
$\bo$, $\tor$. In our further considerations, {\sl we will use the functions 
$\theta_0 (z), $ $Z(z),$ which are defined with the parameter $k=\kappa .$}
In addition, we will need ${\rm sn}  (z;\kappa ),$ ${\rm cn}  (z;\kappa )$, 
${\rm dn}  (z;\kappa )$, as well as ${\rm sn}  (z;x_0^2),$ 
${\rm cn}(z;x_0^2)$, ${\rm dn}  (z;x_0^2)$; the argument $x_0^2$ appears 
because for $k=\kappa ,$ we have $l=x_0^2 .$ 

\medskip

The function $\theta_0(z)$ is entire and has simple zeros at the points
$$z_{m,n}=\iu K'+2mK+2\iu nK',\quad\text{ for }\,\,\,m,n\in\Z ;$$
likewise, $Z(z)$ is a  meromorphic function with the simple poles at the 
points $z_{m,n}$, for $m,n\in\Z$. In addition, the functions $\theta_0$ and 
$Z$ are ``almost'' double-periodic:
\begin{eqnarray*}
\theta_0(z+2K)=\theta_0 (z) ,\quad
\theta_0(z+2\iu K')=-h^{-1}e^{-\frac{\pi \iu z}{K}}\theta_0(z), \\
Z(z+2K)=Z(z),\quad Z(z+2\iu K')=Z(z)-\frac{\pi \iu}{K} .
\end{eqnarray*}

\medskip


We consider the rectangle
$$
\funddomain =\bigl\{ z\in\C : \quad -2L< {\rm Re}\,   z< 2L, \quad -L'<{\rm Im} \,  
z< L'\bigr\},
$$
and the analytic function
$$
\bsigma(z) =x_0 \, {\rm sn}  ( z+L ;x_0^2) .
$$
Let us introduce the following subrectangles of $\funddomain$:
\begin{eqnarray*}
\funddomain^{-}& =&\bigl\{ z\in\C:\quad -2L<{\rm Re} \,   z<0,\ \ 
-L'<{\rm Im} \, z<L'\bigr\},\\
\funddomain^{+}& =&\bigl\{ z\in\C:\quad 0<{\rm Re} \,   z<2L,\ \ 
-L'<{\rm Im} \,  z<L'\bigr\},\\
\funddomain^{-}_1& =&\bigl\{ z\in\C:\quad -2L<{\rm Re} \,   z<0,\ \ 
-L'<{\rm Im} \,  z<0\bigr\}\\
\funddomain^{-}_0&=&\bigl\{ z\in\C:\quad  -2L<{\rm Re} \,   z<0,\ \ 
0<{\rm Im} \,  z<L'\bigr\},\\
\funddomain^{+}_1& =&\bigl\{ z\in\C:\quad 0<{\rm Re} \,   z<2L,\ \ 
-L'<{\rm Im} \,  z<0\bigr\} \\
\funddomain^{+}_0&=&\bigl\{ z\in\C:\quad 0<{\rm Re} \,   z<2L,\ \ 0<{\rm Im} 
\,  z<L'\bigr\} .
\end{eqnarray*}
The function $\bsigma$ maps each of the rectangles $\funddomain^{-}$ and 
$\funddomain^{+}$ conformally onto the slit sphere
$$\s\setminus\Big([-x_0 ;x_0]\,\,\cup\,\,] -\infty ;-1/x_0]
\cup [1/x_0;+\infty [ \,\,\cup\,\{\infty\}\Big).$$
It is also known that $w=\bsigma(z)$ maps the closed rectangle 
$\bar{\funddomain}$ with identified opposite sides conformally onto 
$\tor$ (see \cite[Ch. VIII]{Akhiez} or \cite[Ch. VI, pp. 280--285]{Neehari}).

The inverse function of the restriction of $w=\bsigma(z)$ to 
$\funddomain^{-}$ is given by the elliptic integral
$$
z=\btau (w)=\int\limits_{0}^{w/x_0}
\frac{\diff t}{\sqrt{(1-t^2)(1-x_0^4t^2)}} -L .
$$
As a conformal mapping, $z=\btau (w)$ sends the upper half-plane
$$\C_{+} =\big\{ w\in\C:\,\, {\rm Im} \,  w>0\big\}$$
onto the rectangle $\funddomain_0^-$, and it sends the lower half-plane
$$\C_{-}=\big\{ w\in\C:\,\, {\rm Im} \,  w<0\big\}$$
onto the rectangle $\funddomain_1^-$, in such a way that
\begin{gather*}
\btau (x_0 )=0, \qquad \btau (-x_0) =-2L,\qquad \btau (0) =-L,\\
  \lim_{\C_{+}\ni w\to 1/x_0} \btau  (w) =\iu L',   \qquad
 \lim_{\C_{-}\ni w\to 1/x_0} \btau  (w) =-\iu L' ,
\\
\lim_{\C_{+} \ni w\to -1/x_0} \btau  (w) =-2L+\iu L',  \qquad
\lim_{\C_{-}\ni w\to -1/x_0} \btau  (w) =-2L-\iu L' ,
\\
\lim_{\C_{+} \ni w\to\infty}\btau (w )=-L+\iu L' ,\qquad
\lim_{\C_{-}\ni w\to\infty}\btau (w )=-L-\iu L'.
\end{gather*}
The function $z=\btau (w)$ extends to an analytic function on
$$\C_+\cup\C_-\,\cup\,\,]x_0,1/x_0[\,\,\cup\,\,]-1/x_0,-x_0[.$$
Its restriction to the upper half plane $\C_+$ has an analytic continuation 
across the remaining segments
$$\R\cup\{\infty\}\setminus \Big(]x_0,1/x_0[\,\,\cup\,\,]-1/x_0,-x_0[\Big),$$
and so does its restriction to the lower half plane $\C_-$. If we look 
carefully at these extensions, we find that the mapping $z=\btau (w)$ lifts 
to a conformal mapping $\tor\to\C/\bgamma$, where $\bgamma$ is the additive 
group generated by the elements $4L$ and $2\iu L'$.
We let $\funddomain_{\rm fund}$ denote the set $\funddomain$ adjoined with 
the left vertical and the lower horizontal sides of this rectangle;
then, $\funddomain_{\rm fund}$ is a fundamental domain for $\C/\bgamma .$

We need understand the operations $p\mapsto p'$ and $p\mapsto p^*$ on $\tor$ 
in terms of this
identification of $\tor$ with $\C/\bgamma$. It is easy to see that the 
mirror mapping $p\mapsto p'$
corresponds to $z\mapsto -z$ on $\C/\bgamma$. Also, the reflection in 
$\partial\Delt$ mapping
$p\mapsto p^*$ corresponds to $z\mapsto z^*$, where $z^*$ is the reflected 
point in the line
$\frac{\iu}2L'+\R$ (modulo $\bgamma$). This latter fact is perhaps not 
entirely obvious.
To see that it is nevertheless so, pick a point $z\in\funddomain_{\rm fund}$. 
We have
$$
\bsigma(z)= w,\qquad {\rm sn}  (z+L;x_0^2)  =\frac{w}{x_0}.
$$
Using the relation (see \cite[table XII]{Akhiez})
$$
\mbox {sn} (u+\iu L';x_0^2 ) =\frac{1}{x_0^2\,{\rm sn}  (u;x_0^2 )} ,
$$
we find that
$$
{\rm sn}  (\bar{z}+L+\iu L';x_0^2) =\frac{1}{x_0^2{\rm sn}  
(\bar{z}+L;x_0^2)}= \frac{1}{x_0\bar{w}} ,
$$
so that
$$x_0\, {\rm sn}  ((\bar{z}+\iu L')+L;x_0^2) =w^*.$$
We realize that $$z^*=\bar{z}+\iu L' ,$$
which is the formula expressing reflection in the line $\frac{\iu}2L'+\R$.

Finally, we obtain a description of the image $\btau (\Delt ):$
{\sl it is the rectangle }
$$
\funddomain_{\Delt} =\left\{ z\in\C:\ -2L\le {\rm Re} \, z<2L,\,\,\, 
|{\rm Im} \,  z|<\frac{L'}{2}\right\} .
$$

The image of $\partial \Delt$ consists of the two horizontal line segments
$$\gamma_{\pm} =\biggl\{ -2L\le{\rm Re} \,   z< 2L,\ {\rm Im} \, 
 z=\pm\frac{L'}{2}\biggr\}.$$

\medskip
For $ (z,\zeta )\in\C\times \C $, we define the
function $G(z,\zeta)$ by
\begin{equation*}
G (z,\zeta) = -\frac{1}{2}\log\left|
\frac{\theta_0(Mz-M\zeta +\iu K')\theta_0(M\bar{z}+M\zeta)}{\theta_0 
(Mz+M\zeta-\iu K')
\theta_0 (M\bar{z}-M\zeta )}\right|^2
-\frac{\pi M}{K} \left[\frac{2M}{K'}{\rm Im} \,  \zeta -1\right]{\rm Im} 
\,  z ,
\end{equation*}
where the function $\theta_0$ is given by (\ref{el-J}); here, we think of 
$\log$ as taking values in $[-\infty ;+\infty ].$

From the properties of the function  $\theta_0$, and (\ref{L-M}), we can 
easily obtain that $G(z,\zeta)$ has the following properties:

\medskip

\noindent
$1^{\circ}$ $G(z,\zeta )=G(\zeta ,z);$

\medskip

\noindent
$2^{\circ}$ the function $z\mapsto G(z,\zeta )$ is periodic with respect 
to the group $\bgamma ,$ making it a  function on $\C/ \bgamma;$

\medskip

\noindent
$3^{\circ}$ for a fixed $\zeta\in\funddomain_{\Delt },$ the function 
$z\mapsto G(z,\zeta )$ is harmonic in the variable $z$
in the domain $\funddomain_{\Delt}\setminus \{ \zeta, -\zeta \} ,$
it has the logarithmic singularities $\log |z-\zeta |+O(1)$ near $z=\zeta $
and $-\log |z+\zeta |+O(1)$ near $z=-\zeta$;

\medskip

\noindent
$4^{\circ}$ $G(z,\zeta )=0$ as $z\in\gamma_{+}\cup\gamma_- ;$

\medskip

\noindent
$5^{\circ}$ $G(-z,\zeta )=-G(z,\zeta ).$

\bigskip

The property $2^{\circ}$ means that $G(\btau (p), \btau (q))$ is a function 
on $\tor\times\tor .$
From the above properties of $G$, it also follows that $G(\btau  (p) ,
\btau  (q))$ coincides
with the previously considered function $G_{\Delt}^{{\rm alt}} (p,q)$:
\begin{equation}
G_{\Delt}^{\rm alt} (\bsigma(z),\bsigma(\zeta)) \equiv G(z,\zeta ),\qquad
(z,\zeta )\in\C/ \bgamma\times\C/ \bgamma .
\label{G_Delt}
\end{equation}

We denote by $\bpi_{\Delt}$ the subdomain of the torus $\C/\bgamma$ whose
restriction to the fundamental domain $\funddomain_{\rm fund}$ is the 
subrectangle $\funddomain_{\Delt}$, and by $G_{\bpi_{\Delt}} (z,\zeta)$ 
the Green function of this subdomain. Then, the relation (\ref{G_Delt}) 
is equivalent to
$$
G(z,\zeta )=G_{\bpi_{\Delt}}^{\rm alt}(z,\zeta )=G_{\bpi_{\Delt}} (z,\zeta )
-G_{\bpi_{\Delt}} (-z,\zeta ),\qquad (z,\zeta )\in
\bpi_{\Delt}\times\bpi_{\Delt} .
$$
Let us consider the function
\begin{equation}
Q_{\Delt} (z) =\partial_{\zeta} G(z,\zeta)\vert_{\zeta =0}
=MZ(Mz+\iu K')-MZ(M\bar{z}) +\frac{\pi \iu M}{KK'}{\rm Im} \,  (Mz) + 
\frac{\pi \iu M}{2K},
\label{def-Q}
\end{equation}
where $Z$ is Jacobi  $Z$-function (see (\ref{el-J})).
The above properties of $G(z,\zeta )$ imply that $Q_{\Delt} (z)$ has
the properties:

\medskip

\noindent
(1) it is periodic function with respect to $\bgamma$, so that 
$Q_{\Delt} (z)$ is a function on $\C /\bgamma ;$

\medskip

\noindent
(2)  the function $Q_{\Delt}$ is harmonic on 
$\funddomain_{\Delt}\setminus\{ 0\}$, and it has the singularity 
$1/z+O(1)$ at the point $0$;

\medskip

\noindent
(3) $Q_{\Delt}(z)=0$ for $z\in \gamma_{+}\cup\gamma_- ;$

\medskip

\noindent
(4) $Q_{\Delt}(-z)=-Q_{\Delt}(z) ,$ for  $ z\in\C .$

\bigskip

Put
$$
Q_1 (\bp ) \equiv \big( Q_{\Delt}\circ \btau  \circ\bphi\big) (\bp) ,
\qquad \bp\in\bo .
$$
This function satisfies the conditions (Q1) and (Q2) of 
Proposition~\ref{prop2-3}.
Also, it has the singularity $$\frac{1}{\btau (\phi (z^2))}+O(1)\sim 
\frac{b}{z} +O(1)$$
at the point $\bzero \in\bo ;$
here,
\begin{multline}
b=\lim\limits_{z\to 0}\frac{z}{\btau (\phi (z^2))} =
\lim\limits_{w\to x_0}\frac{\sqrt{\varphi (w)}}{\btau (w)}\\
=\lim\limits_{w\to x_0}\frac{\sqrt{w-x_0}}{\btau (w)} =
-\lim\limits_{w\to x_0}\frac{(\sqrt{w-x_0})'_w}{ \btau '(w)} 
=\frac{\iu}{\sqrt{2}}\sqrt{x_0 (1-x_0^4)} .
\label{b}
\end{multline}
In view of the above, it follows that
$$
Q (\bp ) =\frac{1}{b}\, Q_1(\bp ) =\frac{1}{b}\, 
\big( Q_{\Delt}\circ \btau  \circ\bphi\big) (\bp)
$$
is exactly the function we are looking for.
\bigskip

\noindent
{\bf The area-theorem type inequality for univalent function on $\D$.}
We now write down the inequality (\ref{2-3}) for the function
$$P(z)=(R \circ \bvarphi \circ \bsigma)(z) -
\frac{1}{b}\,Q_{\Delt} (z),\qquad z\in\funddomain_{\Delt} .$$
As we recall the definition of the function $R,$ we see that
$$(R \circ \bvarphi \circ \bsigma)(z) =\frac{1}{\sqrt{\varphi 
(\bsigma(z))}},$$
where $\sqrt{u}$ means the algebraic square root of $u.$
Then, for our choice of $P$, the inequality (\ref{2-3})  assumes the form
\begin{eqnarray}
\int\limits_{\funddomain_{\Delt}}\left|
\frac{\varphi'(\bsigma(z))\bsigma'(z)}{2\big[\varphi 
(\bsigma(z))\big]^{3/2}} + \frac{1}{b}\,\partial_zQ_{\Delt} (z)
\right|^2 \, \diff A (z)
 \le \frac{1}{|b|^2}
\int\limits_{\funddomain_{\Delt}} \big|\bar\partial_z Q_{\Delt} (z)\big|^2
\, \diff A (z);
\label{3-6}
\end{eqnarray}
here, as usual, $\diff A (z)$ is the area element, and the constant $b$ 
is as in (\ref{b}).

\medskip

We intend to simplify the inequality (\ref{3-6}).
First, we evaluate the right-hand side of (\ref{3-6}).
Let us recall that
$$
Q_{\Delt} (z) =\partial_{\zeta} G(z,\zeta )\big\vert_{\zeta =0} =
\partial_{\zeta} \bigg\{ G_{\bpi_{\Delt}} (z,\zeta) 
-G_{\bpi_{\Delt}} (-z,\zeta)\bigg\}\bigg\vert_{\zeta =0},
$$
so that
$$
\bar\partial_zQ_{\Delt} (z) = \bigg\{\bar\partial_z
\partial_{\zeta} G_{\bpi_{\Delt}} (z,\zeta)+
\bar\partial_z\partial_{\zeta}G_{\bpi_{\Delt}} (-z,\zeta)\bigg\}
\bigg\vert_{\zeta =0} .
$$

The kernel
$$K_{\bpi_{\Delt}}(z,\zeta )=-\frac{2}{\pi}\, 
\partial_z\bar\partial_{\zeta}\, G_{\bpi_{\Delt}}(z,\zeta ),
\qquad z\neq\zeta ,$$
has the following reproducing property: {\sl for any analytic function 
$f\in L^2(\bpi_{\Delt})$},
$$f(\zeta ) =\int\limits_{\funddomain_{\Delt}} f(z)
\bar K_{\bpi_{\Delt}}(z,\zeta )\, \diff A(z),
\qquad \zeta\in\bpi_{\Delt}.$$
In particular, taking into account that the function  
$z\mapsto K_{\bpi_{\Delt}} (z,\zeta )$ is analytic and bounded near
the point $z=\zeta$, we have
$$
\int\limits_{\funddomain_{\Delt}} |K_{\bpi_{\Delt}}(z,\zeta )|^2 \diff A(z) = 
K_{\bpi_{\Delt}}(\zeta ,\zeta ).
$$
From the above, it follows that
\begin{multline}
\int\limits_{\funddomain_{\Delt}} \big|\bar\partial_z Q_{\Delt} (z)\big|^2
\, \diff A (z) =\frac{\pi^2}{4} \int\limits_{\funddomain_{\Delt}} 
\big|\bar K_{\bpi_{\Delt}}(z,0)+
\bar K_{\bpi_{\Delt}}(-z,0)\big|^2\diff A(z) \\
=\pi^2K_{\bpi_{\Delt}}(0,0) =
-2\pi\,  \partial_z\bar\partial_{\zeta}\, G_{\bpi_{\Delt}}(0,0)=
 -\pi\, \overline{\bar\partial_z\,Q_{\Delt } (0) }= -\pi\, 
\partial_z\,\bar Q_{\Delt } (0).
\label{auxil}
\end{multline}

The calculations of the right-hand side of (\ref{3-6}) can be completed by
using the following facts from the elliptic functions theory 
(see \cite[Ch.V]{Akhiez}):
\begin{eqnarray}
Z'(u) =\big[{\rm dn}\, (u;\kappa )\big]^2-\frac{E}{K},
\label{3-7}
\\
{\rm dn}\,  (0;\kappa ) =1,\nonumber \\
EK'+E'K-KK'=\frac{\pi}{2},
\label{3-9}
\end{eqnarray}
where $E=E(\kappa)$, $E'=E(\kappa')$, $K=K(\kappa)$, and $K'=K(\kappa')$
(see equations (\ref{el-E}), (\ref{el-K}), and (\ref{M})).  
In view of (\ref{def-Q}), we have
\begin{equation}
\bar\partial_zQ_{\Delt } (z) =M^2\left( -\big[ {\rm dn}\, 
( M\bar{z};\kappa )\big]^2 +\frac{E}{K}
-\frac{\pi}{2KK'}\right) ,
\label{der-Q}
\end{equation}
so that
\begin{equation*}
\bar\partial_zQ_{\Delt } (0) =-\frac{M^2E'}{K'} ,
\end{equation*}
which is a real number.
We get, by (\ref{auxil}),
\begin{equation}
\frac{1}{|b|^2}\int\limits_{\funddomain_{\Delt}} 
\big|\bar\partial_z Q_{\Delt} (z)\big|^2
\,\diff A (z) =\frac{\pi M^2E'}{|b|^2K'}
=\frac{\pi (1+x_0^2)E'}{2x_0(1-x_0^2)K'} .
\label{point0}
\end{equation}
Finally, the inequality (\ref{3-6}) becomes
\begin{equation}
\int\limits_{\funddomain_{\Delt}}\left|
-\frac{\varphi'(\bsigma(z))\bsigma'(z)}
{2\big[\varphi (\bsigma(z))\big]^{3/2}} -\frac{1}{b}\partial_zQ_{\Delt} (z)
\right|^2 \, \diff A (z)
 \le\frac{\pi M^2E'}{|b|^2K'}= \frac{\pi (1+x_0^2)E'}{2x_0(1-x_0^2)K'} .
\label{3-6new}
\end{equation}
\bigskip

\noindent
{\bf A pointwise estimate.} Put
$$
\bbpsi (z) =-\frac{\varphi'(\bsigma(z))\bsigma'(z)}
{2\big[\varphi (\bsigma(z))\big]^{3/2}}
-\frac{1}{b}\,\partial_zQ_{\Delt} (z) ,\qquad z\in\bpi_{\Delt};
$$
the inequality (\ref{3-6new}) now takes form
\begin{equation}
\int\limits_{\funddomain_{\Delt}} |\bbpsi (z)|^2 \, \diff A (z)\le 
\frac{\pi M^2 E'}{|b|^2 K'} .
\label{3-10}
\end{equation}
By (\ref{point0}), (\ref{3-10}), and the reproducing 
property of the function
$$-\frac{1}{\pi }\overline{\bar\partial_{z}
Q_{\Delt} (z) }=\frac{1}{2} \big(K_{\bpi_{\Delt}}(z,0) +
K_{\bpi_{\Delt}}(-z,0)\big),$$
we have -- by the Cauchy-Schwarz inequality --
\begin{multline*}
|\bbpsi (0)|^2 =\bigg| \, \int\limits_{\funddomain_{\Delt}} \bbpsi (z)
\left[-\frac{1}{\pi}\,\bar\partial_{z}Q_{\Delt} (z) \right]
\, \diff A (z) \bigg|^2\\
\le\frac1{\pi^2}\int\limits_{\funddomain_{\Delt}} | \bbpsi (z)|^2\diff A (z)
\int\limits_{\funddomain_{\Delt}} \big|\partial_{z}\bar Q_{\Delt} (z)\big|^2
\diff A (z)\le \frac{M^4}{|b|^2} \left( \frac{E'}{K'}\right)^2,
\end{multline*}
whence,
\begin{equation}
|\bbpsi (0)| \le \frac{M^2}{|b|}\cdot\frac{E'}{K'} .
\label{3-11}
\end{equation}
Below, we shall demonstrate that this inequality is equivalent to the  
estimate (\ref{eq-ptests}) of  Goluzin for the class $\Sigma$.

\bigskip

\noindent
{\bf Rewriting the area-type inequality in the coordinates of the unit disk.}
We first rewrite the left-hand side of (\ref{3-6new}) as an integral over 
$\Delt$ rather than $\funddomain_{\Delt}$:
\begin{multline}
\int\limits_{\funddomain_{\Delt}}\bigg|
-\frac{\varphi'(\bsigma(z))\bsigma'(z)}{2\big[\varphi 
(\bsigma(z))\big]^{3/2}} -\frac{1}{b}\partial_zQ_{\Delt} (z)
\bigg|^2\,  \diff A (z) =  \\
\int\limits_{\Delt}\bigg|
-\frac{\varphi' (w)}{2\big[\varphi (w)\big]^{3/2}} -\frac{1}{b}\, 
\partial_zQ_{\Delt} (z)\Bigr \vert_{z=\btau (w)} \,\btau '(w)
\bigg|^2 \, \diff A (w);
\label{3-12}
\end{multline}
here, the area measure $\diff A$ is implicitly lifted  from $\D$ to $\Delt$.
From (\ref{3-7}), (\ref{3-9}) and the following relations between Jacobi 
elliptic functions (\cite[table XII]{Akhiez}),
\begin{eqnarray*}
{\rm dn}\, (u+\iu K';\kappa ) =-\iu\frac{{\rm cn}\, (u;\kappa )}{{\rm sn}\, 
( u;\kappa )} ,\\
\big[{\rm sn}\, (u;\kappa )\big]^2+\big[ {\rm cn}\, (u;\kappa )\big]^2=1,
\end{eqnarray*}
we obtain, in view of  (\ref{def-Q}),
$$
\partial_zQ_{\Delt} (z) =M^2\left[ -\frac{1}{\big[{\rm sn}\, 
(Mz;\kappa )\big]^2} +\frac{E'}{K'}\right] .
$$
We note that the expression
$$\big[{\rm sn}\,  ( Mz;\kappa)\big]^2\big\vert_{z=\btau (w)}$$
can be simplified by using Landen's transformation of Jacobi functions
(\cite[Ch.VI]{Akhiez}).
This transformation  allows us to express  $\big[{\rm sn}\, 
(Mz;\kappa)\big]^2$ as
a function of the expression
 $$\xi (z) =\frac{{\rm cn}\, (z;x_0^2)}{{\rm dn}\, (z;x_0^2)}.$$
We have
\begin{eqnarray*}
{\rm cn}\, (z;x_0^2) =\frac{1-(1+\kappa')\big[{\rm sn}\, 
(Mz;\kappa)\big]^2}{{\rm dn}\, (Mz;\kappa)},\\
{\rm dn}\, (z;x_0^2) =\frac{1-(1-\kappa')\big[{\rm sn}\, 
(Mz;\kappa)\big]^2}{{\rm dn}\, (Mz;\kappa)} .
\end{eqnarray*}
From these formulas we find that
$$
\big[{\rm sn}\, (Mz;\kappa)\big]^2=
\frac{1-\xi (z)}{1+\kappa'-(1-\kappa')\xi (z)} .
$$
Further, taking into account the relation
$${\rm sn}\, (z+L;x_0^2) =
\frac{{\rm cn}\,  (z;x_0^2)}{{\rm dn}\, (z;x_0^2)},$$
we conclude that $w=\xi (z)$ is the inverse function to 
$$w\mapsto \int\limits_{0}^{w}\frac{\diff t}{\sqrt{(1-t^2)(1-x_0^4t^2)}}.$$
From the above, we obtain that
$$
\big[{\rm sn}\, (Mz;\kappa)\big]^2\big\vert_{z=\btau (w)} 
=\frac{1-w/x_0}{1+\kappa'-(1-\kappa')w/x_0} =
\frac{(1+x_0^2)(w-x_0)}{2x_0 (x_0 w-1)}.
$$
Finally, we arrive at
$$
\partial_zQ_{\Delt} (z)\Bigr\vert_{z=\btau (w)} =
\frac{x_0 (1+x_0^2)}{2} \frac{1-x_0 w}{w-x_0} +\frac{(1+x_0^2)^2E'}{4K'} ,
\qquad w\in\Delt .
$$
Note that the above expression is a well-defined function on $\D$.
Substituting the last expression as well as 
$$\btau'(w)=\frac{1}{\iu\sqrt{(w^2-x_0^2)(1-x_0^2w^2)}}$$
into the right-hand side integral of  (\ref{3-12}),
we get
\begin{multline*}
 \int\limits_{\Delt}
\bigg|-
\frac{\varphi' (w)}{2\big[\varphi (w)\big]^{3/2}} -\frac{1}{b} 
\partial_zQ_{\Delt} (z)\Bigr \vert_{z=\btau (w)} \, \btau '(w)
\bigg|^2 \, \diff A (w)\\
=\int\limits_{\Delt}\biggl|
\frac{\varphi' (w)}{2\big[\varphi (w)\big]^{3/2}} -
\frac{(1+x_0^2)\sqrt{x_0}}{\sqrt{2(1-x_0^4)}}\sqrt{\frac{1-x_0 w}{1+x_0 w}}
\frac{1}{\sqrt{w+x_0 }}\frac{1}{(w-x_0)^{3/2}} \\
 -\frac{(1+x_0^2)^2}{2\sqrt{2x_0(1-x_0^4)}}\frac{E'}{K'}
\frac{1}{\sqrt{(w+x_0)(1-x_0^2w^2) }}\frac{1}{(w-x_0)^{1/2}}
\biggr| ^2\, \diff A (w) .
\end{multline*}
As we multiply by $\sqrt{w^2-x_0^2}$ inside the absolute value signs of 
the integral and divide by $|w^2-x_0^2|$ outside them, which permits us to 
integrate over $\D$ instead of over the covering surface $\Delt$, we 
realize that we have derived the following from (\ref{3-6new}).

\begin{prop}
Let  $\varphi :\D\to\s$ be a univalent function with the following property:
for some real $x_0$, $0<x_0<1,$ we have $\varphi (x_0 )=0$, 
$\varphi (-x_0 )=\infty$, and $\varphi' (x_0 )=1$. Then
\begin{multline}
\int\limits_{\D}
\Biggl|
\frac{\varphi' (w)\sqrt{w^2-x_0^2}}{\big[\varphi (w)\big]^{3/2}} -
\frac{(1+x_0^2)\sqrt{2x_0}}{\sqrt{(1-x_0^4)}}\, \sqrt{\frac{1-x_0 w}{1+x_0 w}}
\, \frac{1}{w-x_0}\\
 -\frac{E'}{K'}\, \frac{(1+x_0^2)^2}{\sqrt{2x_0(1-x_0^4)}}\,
\frac{1}{\sqrt{1-x_0^2w^2 }}
\Biggr| ^2\frac{ \diff A (w) }{|w^2-x_0^2|}\le
\frac{\pi E'}{K'}\, \frac{1+x_0^2}{x_0(1-x_0^2)} ,
\label{3-13}
\end{multline}
where $E'=E\bigl( (1-x_0^2)/(1+x_0^2)\bigr) ,$ $K'=
K\bigl( (1-x_0^2)/(1+x_0^2)\bigr)$ and the functions $E(\lambda )$,
$K(\lambda )$ are defined by $(\ref{el-E})$ and $(\ref{el-K})$.
Equality is attained in $(\ref{3-13})$ if and only if  $\varphi$ is a 
full mapping.
\end{prop}

\medskip

\noindent
{\bf The corresponding estimates the class $\Sigma$.}
Let $\psi (z) =z+b_0 +b_1z^{-1} +\dots$ be an element of the class $\Sigma .$
Fix a point $\zeta\in\D_e\setminus\{\infty\} .$
Then
\begin{equation}
x_0 =\frac{1-\sqrt{1-|\zeta |^{-2}}}{1+\sqrt{1-|\zeta |^{-2}}}
\label{x-0}
\end{equation}
satisfies $0<x_0<1$ and we have the inverse relation
$$|\zeta | =\frac{1+x_0^2}{2x_0}.$$
The mapping  
$$\eta (z)=\frac{|\zeta|-x_0\bar\zeta z}{\bar\zeta z-x_0|\zeta|}$$ 
maps $\D_e$ onto $\D$ conformally and takes $\infty$ to $-x_0$  while $\zeta$ 
is mapped to $x_0$. The inverse mapping is 
$$\eta^{-1}(w)=\frac{\zeta}{|\zeta|}\,\frac{1+x_0w}{w+x_0}.$$
Consider the related function
\begin{eqnarray}
\varphi (w) =\frac{|\zeta|}{\zeta}\,\frac{(1+x_0^2)^2}{1-x_0^2}  
\,\frac{\psi(\eta^{-1}(w))-\psi (\zeta )}{\zeta^2\,\psi'(\zeta )} ,
\label{3-14}
\end{eqnarray}
which is  univalent  on $\D$ with $\varphi (-x_0 )=\infty ,$
$\varphi (x_0 )=0,$ $\varphi' (x_0 )=1 .$

Substituting  (\ref{3-14}) into (\ref{3-13}) and making the change of 
variable $w=\eta (z),$  we obtain, after some  simplification,
the corresponding inequality for $\psi$. We write it down in the following 
form.

\begin{thm}{\rm (The area-type estimate)}
Fix a point $\zeta\in\D_e\setminus\{\infty\}$.
Then, for any $\psi\in\Sigma$,
\begin{multline}
\int\limits_{\D_e} \biggl|
\left( \frac{\psi' (\zeta )(z-\zeta )}{\psi (z) -\psi (\zeta )}\right)^{1/2}
\frac{\psi' (z)}{\psi (z)-\psi (\zeta )} -
\left( \frac{1-(\bar\zeta z)^{-1}}{1-|\zeta|^{-2}}\right)^{1/2}   
\frac{1}{z-\zeta}\, \\
+\frac{E'}{K'}\, 
\frac{1}{\big[ (1-|\zeta|^{-2})(1-(\bar\zeta z)^{-1})\big]^{1/2}\,z}
\biggr|^2\, \frac{\diff A(z)}{|z-\zeta |}
\le\frac{2\pi E'}{K'}\, \frac{|\zeta |}{|\zeta |^2-1} ,
\label{3-15}
\end{multline}
where $E'=E\bigl(\sqrt{1-|\zeta |^{-2}}\bigr)$, 
$K'=K\bigl(\sqrt{1-|\zeta |^{-2}}\bigr)$ and the functions $E(\lambda)$,
$K(\lambda)$ are defined by  $(\ref{el-E})$ and $(\ref{el-K})$.
The above inequality is an equality if and only if $\psi$ is a full 
mapping.
\end{thm}

\medskip

\noindent
{\bf The derivation of Goluzin's inequality from the area-type estimate.}
Put
\begin{multline*}
\Psi (z,\zeta )=
\left( \frac{\psi' (\zeta )(z-\zeta )}{\psi (z) -\psi (\zeta )}\right)^{1/2}
\frac{\psi' (z)}{\psi (z)-\psi (\zeta )} -
\left(\frac{1-(\bar\zeta z)^{-1}}{1-|\zeta|^{-2}}\right)^{1/2}   
\frac{1}{z-\zeta} \\
+\frac{E'}{K'} \frac{1}{\big[(1-|\zeta|^{-2})(1-(\bar\zeta z)^{-1})
\big]^{1/2}\,z}.\end{multline*}
As we recall how the inequality (\ref{3-10}) containing the function 
$\bbpsi$ is transformed into (\ref{3-15}) involving the analogous function
$\Psi$, we find that
$$
|\Psi (\zeta ,\zeta ) |=\frac{x_0^{3/2}\sqrt{1-x_0^4}}{\sqrt{2}}\,
\frac{|\zeta|^2}{|\zeta|^2-1}\, |\bbpsi (0)|.
$$
In view of (\ref{3-11}), we then have
\begin{equation*}
|\Psi  (\zeta ,\zeta )|\le \frac{M^2E'}{\sqrt{2}\,|b|\,K'} \,
x_0^{3/2}\sqrt{1-x_0^4}\,\frac{|\zeta|^2}{|\zeta|^2-1} ,
\end{equation*}
where $x_0$ is given in terms of $|\zeta|$ by (\ref{x-0}).
By substituting the expressions for the constants $M$ and $b$ 
(see (\ref{M}) and (\ref{b})), and simplifying further,
we obtain the estimate
\begin{equation}
|\Psi  (\zeta ,\zeta )|\le\frac{E'}{K'}\,\frac{|\zeta|}{|\zeta |^2-1} .
\label{3-16}
\end{equation}
On the other hand, a direct calculation yields
$$
\Psi (\zeta, \zeta ) =
\frac{\psi'' (\zeta )}{4\psi' (\zeta )}
-\frac{1}{2\zeta}-\frac{2-|\zeta |^2}{2\,(|\zeta |^2-1)\,\zeta}
+\frac{E'}{K'}\,\frac{|\zeta|^2}{(|\zeta |^2-1)\,\zeta } .
$$
The inequality (\ref{3-16}) thus takes the form
\begin{equation}
\left|\frac{\zeta\,\psi'' (\zeta )}{\psi' (\zeta )}
-2+\frac{2(|\zeta |^2-2)}{|\zeta |^2-1}
+\frac{4E'}{K'}\frac{|\zeta|^2}{ (|\zeta |^2-1)}\right|\le
\frac{E'}{K'}\frac{4|\zeta|^2}{|\zeta|^2-1} .
\label{3-17}
\end{equation}
From (\ref{3-9}), we have
$$
\frac{E'}{K'} =1 -\frac{E}{K} +\frac{\pi }{2KK'} ,
$$
which, together with (\ref{3-17}), leads to
$$
\left|
\frac{\zeta\, \psi''(\zeta )}{\psi'(\zeta )} +
\frac{4|\zeta |^2 -2}{|\zeta |^2 -1} -\frac{4|\zeta |^2}{|\zeta |^2-1 }
\frac{E\left(\frac{1}{|\zeta |}\right)}{K\left(\frac{1}{|\zeta |}\right)}
\right| \le \frac{4|\zeta |^2}{|\zeta |^2-1}
 \left( 1-\frac{E\left(\frac{1}{|\zeta |}\right)}{K 
\left(\frac{1}{|\zeta |}\right)}\right)
$$
after some simplification;
here, the functions $E(\lambda ),$ $K(\lambda )$ are defined by  
(\ref{el-E}) and (\ref{el-K}).
This is the classical inequality due to Goluzin 
(see \cite{Gol1943}, \cite[Ch.IV, \S 3, p. 132]{Goluz}), and if we 
divide by $\zeta$ inside the absolute value parentheses, we arrive at 
(\ref{eq-ptests}).

\begin{rem}
To find the extremal $\psi\in\Sigma$ which gives equality in Goluzin's 
inequality (\ref{eq-ptests}) at a given point $z\in\D_e$, we should just 
check when we have equality in the Cauchy-Schwarz inequality leading up 
to (\ref{3-11}). The result of this exercise of course agrees with Goluzin's 
findings. 
\end{rem}

\end{document}